\documentclass[12pt,draft,amsmath,amscd,amsbsy,amssymb]{amsart}

\usepackage[english]{babel}

 \relax

\chardef\bslash=`\\

\makeatletter \def\verbatim{\interlinepenalty\@M \@verbatim
\leftskip\@totalleftmargin\advance\leftskip2pc
\frenchspacing\@vobeyspaces \@xverbatim} \makeatother \hfuzz1pc

\makeatletter \def\dgt@k{\dg@DX=-3 \dg@DY=2 \dg@SIZE=3}
\makeatother

\makeatletter \def\dgt@kk{\dg@DX=3 \dg@DY=-1 \dg@SIZE=3}

\theoremstyle{plain}
\newtheorem{thm}{Theorem}[section]

\newtheorem{lemma}[thm]{Lemma}

\theoremstyle{definition}
\newtheorem{defin}[thm]{Definition}

\newcommand{\pr}{{\mathrm p}{\mathrm r}}

\newcommand{\comp}{{\mathrm C}{\mathrm o}{\mathrm m}{\mathrm p} }
\newcommand{\card}{{\mathrm C}{\mathrm a}{\mathrm r}{\mathrm d} }

\newcommand{\G}{{\mathcal G}}
\newcommand{\VG}{{\mathcal V}{\mathcal G}}
\newcommand{\EG}{{\mathcal E}{\mathcal G}}
\newcommand{\Oo}{\mathcal O}

\numberwithin{equation}{section}

\usepackage[all]{xy}

\begin{document}

\title[]
{$\infty$-Open-multicommutativity in the category $\comp$}
\author{Roman Kozhan}

\address{Department of Mechanics and Mathematics, Lviv National University,
Universytetska 1, 79000 Lviv, Ukraine} \email{}

\thanks{The author gratefully thanks Michael Zarichnyi for his capable
assistance in this research.}


\maketitle

\begin{abstract}In this paper the notion of
$\infty$-open-multicommutativity of functors in the category of
compact Hausdorff spaces is considered. This property is a
generalization of the open-multicommutati-vity on the case of
infinite diagrams. It is proved that every open-multicommutative
functor is $\infty$-open-multicommutative.
\end{abstract}

\section*{Introduction}\label{s:intro}

We continue to study the properties of openness and
bicommutativity of functors in the category $\comp$ of compact
Hausdorff spaces. The functor in the category $\comp$ is open if
it preserves the class of open and surjective maps. The question
whether the notions of openness and bicommutativity are equivalent
for normal functors has been set by Shchepin (1981) till now it is
solved only in the case of finite power functor by Zarichnyi
(1992). But in general it is still open.

Due to this problem Kozhan and Zarichnyi (2004) introduce the
notion of open-multicommutativity of normal functors. This
property is a generalization of the bicommutativity extended from
square to more complicated finite diagrams. Simultaneously this
notion includes also the property of openness of the functor which
is necessary condition of the open-multicommutativity. Kozhan and
Zarichnyi (2004) have proved that the functor of probability
measures is open-multicommutative.

There is a natural generalization of this property on infinite
diagrams. Construction of $\infty$-open-multicommutativity, which
is introduced in this note, is the analogous to
open-multicommutativity for the class of infinite diagram. The
main contribution of the paper if the criterium of the
$\infty$-open-multicommutativity for functors in category $\comp$.

In Section \ref{s:defs} we give necessary definitions and notions
concerning to openness and bicommutativity of functors. Section
\ref{s:result} provides the main result of the paper.

\section{Definitions and Preliminaries}\label{s:defs}

Suppose that $\G$ is a finite partially ordered set and we also
regard it as a finite directed graph. Denote by $\VG$ the class of
all vertices of graph $\G$ and by $\EG$ the set of its edges. A
map $\Oo\colon\G\to\comp$ is called a {\em diagram}.

\begin{defin}\label{d:cone} The set of morphisms
\begin{equation}\label{e:1}
(X\overset{g_A}{\to}\Oo(A))_{A \in \VG}
\end{equation}
is said to be a {\em cone} over the diagram $\Oo$ if and only if
for every objects $A,B \in\VG$ and for every edge $\varphi\colon A
\to B$ in $\EG$ the diagram
\[\xymatrix{&X\ar[dl]_{g_A}\ar[dr]^{g_B}&\\\Oo(A)
\ar[rr]_{\Oo(\varphi)}&&\Oo(B)}
\]is commutative.
\end{defin}

\begin{defin}\label{d:limit} The cone (\ref{e:1}) is called a {\em limit} of the
diagram $\Oo$ if the following condition is satisfied: for each
cone $C^\prime=(X^\prime
\overset{{g^\prime_A}}{\to}\Oo(A))_{A\in\VG }$ there exists a
unique morphism $\chi_{C^\prime}\colon X^\prime \to X$ such that
$g^\prime_A=g_A \circ \chi_{C^\prime}$ for every $A\in\VG$.
\end{defin}

Further we denote this cone by $\lim(\Oo)$. The map
$\chi_{C^\prime}$ is called the {\em characteristic map} of
$C^\prime$.

\begin{defin}\label{d:coneom} The cone $C^\prime=(X^\prime
\overset{{g^\prime_A}}{\to}\Oo(A))_{A\in\VG}$ is called {\em
open-multicommutative} if the characteristic map $\chi_{C^\prime}$
is open and surjective.
\end{defin}

Let $F$ be a covariant functor in the category $\comp$. Define the
diagram $F(\Oo)\colon\G\to\comp$ in the following way: for every
$A \in\VG$ let $F(\Oo)(A)=F(\Oo(A))$ and for every edge
$\varphi\in\EG$ we set $F(\Oo)(\varphi)=F(\Oo(\varphi))$.

\begin{defin}\label{d:functom} The functor $F$ is called {\em
open-multicommutative} if it preserves the open-multicommutative
cones, i.e. the cone
\[F(C^\prime)=(F(X^\prime)\overset{F g^\prime_A}{\to}F(\Oo(A)))_{A\in\VG}
\]over the diagram $F(\Oo)$ is open-multicommutative.
\end{defin}

\section{$\infty$-Open-Multicommutativity}
\label{s:result}

We assume, that the graph $\G$ is infinite, i.e. the set $\VG$ is
infinite. Let us denote by
$(Y\overset{\pr_A}{\to}F(\Oo(A)))_{A\in\VG}$ the limit of the
diagram $F(\Oo)$. By the definition \ref{d:coneom} given the cone
$(F(X)\overset{F g_A}{\to} F(\Oo(A)))_{A\in\VG}$ there exists the
unique characteristic map $\chi\colon F(X)\to Y$.

Let $E$ be a set of all finite subsets of $\VG$. Let us define for
every $\alpha \in E$ a finite graph $\G_\alpha$ in the following
way: $\VG_\alpha=\alpha$ and $\varphi$, which connect vertices $A$
and $B$, belongs to $\EG_{\alpha}$ if and only if
$A,B\in\VG_{\alpha}$ ³ $\varphi\in\EG$. Let for every $\alpha \in
E$ the cone $(X_\alpha\overset{g^\alpha_A}{\to}
\Oo(A))_{A\in\VG_\alpha}$ be the limit of the diagram $\Oo$ over
the graph $\G_\alpha$. Define a partial order relation on the set
$E$ as $\alpha\leq\beta\Leftrightarrow\alpha\subseteq\beta$. We
denote by $S=\{X_\alpha,\pr^\beta_\alpha, E\}$ the inverse system.

\begin{lemma}\label{l:lim} $X=\lim S$
\end{lemma}

\begin{proof} Let us show, that the set of morphisms $(X
\overset{\pr_\alpha}{\to}X_\alpha)_{\alpha\in E}$ is a limit of
the inverse system $S$. For arbitrary point $(x_A)_{A\in\VG}\in
X\subseteq\underset{A\in\VG}{\prod}X_A$ its image
$\pr_\alpha((x_A)_{A\in\VG})=(x_A)_{A\in \alpha}\in X_\alpha$ for
each $\alpha\in E$. Indeed, for all $A,B\in \alpha$ and
$\varphi\in\VG_\alpha(A,B)$ we have that
$\pr_B=\Oo(\varphi)\circ\pr_A$ because $A,B\in\VG$ and $\varphi
\in\VG(A,B)$. This implies that $(X\overset{\pr_\alpha}{\to}\Oo(A)
)_{A\in\alpha}$ is a cone. Let us show that this cone is a limit
of $S$. For every cone
$(X^\prime\overset{h_\alpha}{\to}\Oo(A))_{A\in\alpha}$ we define a
map $\zeta\colon X^\prime\to X$ such that
$\zeta(x^\prime)=(h_{\{A\}}(x^\prime))_{a\in\VG}$. Since
$h_\alpha(x^\prime)\in X_\alpha\subseteq\underset{A\in
\alpha}{\prod}\Oo(A)$ for all $\alpha \in E$ we obtain
$h_\alpha(x^\prime) =(h_{\{A\}}(x^\prime))_{A\in\alpha}$ and
therefore $h_\alpha=\pr_\alpha\circ\zeta$. The uniqueness of the
map $\zeta$ proves that the cone
$(X\overset{\pr_\alpha}{\to}\Oo(A))_{A\in\alpha}$ is the limit of
the inverse system $S$.
\end{proof}

Let us consider the cone
$(F(X_\alpha)\overset{F\pr_A}{\to}F(\Oo(A)))_{A\in\alpha}$ over
the diagram $F(\Oo)$ and finite graph $\G_\alpha$. We denote by
$(Y_\alpha\overset{\pr_A}{\to}F(\Oo(A)))_{A\in\alpha}$ the limit
of this diagram over  $\G_\alpha$, $\alpha\in E$. We construct the
inverse systems
\[S_1=\{F(X_\alpha),F(\pr^\beta_\alpha),E\}
\text{ ³ } S_2=\{Y_\alpha,\pr^\beta_\alpha,E\}.
\]Let $\chi\colon
F(X)\to Y$ be a limit map of the morphisms
$(\chi_\alpha)_{\alpha\in E}$, where $\chi_\alpha\colon
F(X_\alpha)\to Y_\alpha$ is a characteristic map of the cone
\[(F(X_\alpha)\overset{F\pr_A}{\to}F(\Oo(A)))_{A\in\alpha}.
\]Lemma \ref{l:lim} implies that $Y=\lim \{Y_\alpha,
\pr^\beta_\alpha, E\}$. It is followed from the uniqueness of the
characteristic maps $\chi_\alpha$ that $\chi$ is a characteristic
map of the cone $(F(X)\overset{F\pr_A}{\to}F(\Oo(A)))_{A\in\VG}$.

\begin{lemma}\label{l:bicdiag} Let $F$ be a bicommutative functor.
For all $\alpha,\beta\in E$ such that $\alpha\leq\beta$ the
diagram
\begin{equation}\label{e:diaglim}
\xymatrix{F(X_\beta)\ar[r]^{\chi_\beta}\ar[d]_{F\pr^\beta_\alpha}&
Y_\beta\ar[d]^{\pr^\beta_\alpha}\\
F(X_\alpha)\ar[r]_{\chi_\alpha}& Y_\alpha}
\end{equation}is bicommutative.
\end{lemma}

\begin{proof} Without loss of generality we assume that
$\card(\beta)=\card(\alpha)+1$ and let $C\in\beta\setminus
\alpha$.

A sufficient condition for the bicommutativity of the diagram
(\ref{e:diaglim}) is the fact that given $\tau_0\in F(X_\alpha)$ ³
$\mu=(\mu_A)_{A\in\beta}\in Y_\beta$ such that
\begin{equation}\label{e:margcond}\chi_\alpha
(\tau_0)=\pr^\beta_\alpha(\mu)=(\mu_A)_{A\in\alpha}
\end{equation} there exists an element $\tau\in F(X_\beta)$ which
satisfies
\[\chi_\beta(\tau)=\mu\text{ and }F\pr^\beta_\alpha(\tau)=\tau_0.
\]
Consider the diagram
\[\xymatrix{F(X_\beta)\ar[rrr]^{F\pr_C}\ar[d]_{F\pr^\beta_\alpha}&&&
F(C)\ar[d]^{{\mathbf 1}_{F(C)}}\\
F(X_\alpha)\ar[rrr]_{{\mathbf 1}_{F(X_\alpha)}}&&& F(\ast)}
\]It is obviously that this diagram is bicommutative, since there exists an element
$\tau\in F(X_\beta)$ such that $F\pr^\beta_\alpha(\tau)=\tau_0$ ³
$F\pr_C(\tau)=\mu_C$. Due to the equation (\ref{e:margcond}) we
see that
\begin{eqnarray*}
\chi_\beta(\tau)&=& \underset{A\in\beta}{\prod}F\pr_A(\tau)=
\underset{A\in \alpha}{\prod}F\pr_A(\tau)\times
F\pr_C(\tau)\\&=&(\underset{A\in \alpha }{\prod}F\pr_A\circ
F\pr_\alpha)\times F\pr_C=(\mu_A)_{A\in \beta}.
\end{eqnarray*}
\end{proof}

\begin{thm}\label{t:wom} Every open-multicommutative functor
$F$ is $\infty$-open-multicommutative.
\end{thm}

\begin{proof} It is known (see Kozhan and Zarichnyi (2004)) that for all $\alpha\in E$
the map $\chi_\alpha$ is open. Lemma \ref{l:bicdiag}
(open-multicommutativity of $F$ implies its bicommutativity) and
Lemma 2.1 of Zarichnyi (2003) imply that the limit map $\chi$ is
also open. The surjectivity of the characteristic map is followed
from the surjectivity of morphisms $\chi_\alpha$, $\alpha \in E$.
This implies that $F$ is $\infty$-open-multicommutative.
\end{proof}

\section*{References}\label{s:bib}

[1] J. Bergin, {\em On the Continuity of Correspondences on Sets
of Measures With Restricted Marginals}, Economic Theory {\bf 13}
(1999), 471-481.

[2] R. V. Kozhan, M. M. Zarichnyi, {\em Open-multicommutativity of
the Probability Measure Functor}, preprint, (2004).

[3] E. Shchepin, {\em Functors and Uncountable powers of
Compacta}, Uspekhi Mat. Nauk. 36, 3, 3-62, (1981) (in Russian).

[4] M. M. Zarichnyi, {\em Characterization of G-symmetric Power
and Extension of Functors onto Kleisli Categories}, Mat. Zametki,
52, 5, (1992), (in Russian).

[4] M. M. Zarichnyi, {\em Correspondences of Probability Measures
With Restricted Marginals Revisited}, preprint, (2003).

\end{document}